\documentclass[12pt]{article}
\usepackage{epsfig,amsfonts,amssymb,latexsym}

\newtheorem{theorem}{Theorem}[section]

\newtheorem{corollary}[theorem]{Corollary}
\newtheorem{proposition}[theorem]{Proposition}

\newenvironment{proof}
{\par\addvspace{0.3cm}\noindent{\rm Proof. }}
{\nopagebreak\mbox{}\hfill $\Box$\par\addvspace{0.25cm}}

\newcommand{\C}{{\mathbb C}}
\newcommand{\Z}{{\mathbb Z}}

\newcommand{\T}{{\mathbb T}}

\newcommand{\G}{F}

\newcommand{\bqn}{\begin{eqnarray}}
\newcommand{\eqn}{\end{eqnarray}}
\newcommand{\ben}{\begin{equation}}
\newcommand{\een}{\end{equation}}
\newcommand{\ba}{\begin{array}}
\newcommand{\ea}{\end{array}}
\newcommand{\nn}{\nonumber}

\newcommand{\iv}{^{-1}}
\newcommand{\iy}{\infty}

\renewcommand{\Re}{{\rm Re\,}}

\newcommand{\diag}{{\rm diag}\,}

\newcommand{\sumarr}[3]{\ba{ccc}\scriptscriptstyle #1\\[-1ex]
\scriptscriptstyle #2\\[-1ex] \scriptscriptstyle #3\ea}

\begin{document}

\date{}
\title{Some identities for determinants of structured matrices} 
\author{Estelle L. Basor \thanks{Supported in part by the
National Science Foundation DMS-9970879.}\\
Department of Mathematics\\ California Polytechnic State
University\\ San Luis Obispo, CA 93407
\and
Torsten Ehrhardt \thanks{Research was in part done at MSRI and
supported by National Science Foundation DMS-9701755.}\\
Fakult\"at f\"ur Mathematik\\ Technische Universit\"at Chemnitz\\
09107 Chemnitz, Germany}

\maketitle
\begin{abstract}
In this paper we establish several relations between the determinants
of the following structured matrices: Hankel matrices, symmetric
Toeplitz + Hankel matrices and Toeplitz matrices. Using known results
for the asymptotic behavior of Toeplitz determinants, these identities
are used in order to obtain Fisher-Hartwig type results on the asymptotics of
certain skewsymmetric Toeplitz determinants and certain Hankel determinants.
\end{abstract}


\section{Introduction}

In this paper we prove identities that involve the determinants of
several types of structured matrices such as Hankel matrices,
symmetric Toeplitz + Hankel matrices and skewsymmetric Toeplitz
matrices.  After having established these identities we show
how they can be used in order to obtain asymptotic formulas for these
determinants.

Let us first recall the underlying notation. Given a sequence
$\{a_n\}_{n=-\iy}^\iy$
of complex numbers, we associate the formal Fourier series
\bqn
a(t)&=&\sum_{n=-\iy}^\iy a_nt^n,\qquad t\in\T.
\eqn
The $N\times N$ Toeplitz and Hankel matrices 
with the (Fourier) symbol $a$ are defined by
\ben
T_N(a) \;\;=\;\; \left(a_{j-k}\right)_{j,k=0}^{N-1},
\qquad
H_N(a) \;\;=\;\; \left(a_{j+k+1}\right)_{j,k=0}^{N-1}.\label{f.THn}
\een
Usually $a$ represents an $L^1$-function defined on the unit circle,
in which case the numbers $a_n$ are the Fourier coefficients,
\bqn
a_n &=& \frac{1}{2\pi}\int_{-\pi}^{\pi}a(e^{i\theta})
e^{-in\theta}\,d\theta,\qquad n\in\Z.
\eqn
To a given symbol $a$ we associate the symbol $\tilde{a}(t):=a(t\iv)$.
The symbol $a$ is called even (odd) if $\tilde{a}(t)=\pm a(t)$, i.e.,
$a_{-n}=\pm a_{n}$.

For our purposes it is important to define another type of Hankel matrix. 
Given a function $b\in L^1[-1,1]$ with moments defined by
\bqn
b_n &=&
\frac{1}{\pi}\int_{-1}^1b(x)(2x)^{n-1}\,dx,\qquad
n\ge1,
\eqn
the $N\times N$ Hankel matrices generated by the (moment) symbol $b$
are defined by
\bqn
H_N[b] &=& \left(b_{1+j+k}\right)_{j,k=0}^{N-1}.
\eqn
We indicate the difference in the definition by using the notation
$H_N(\cdot)$ and $H_N[\cdot]$.
The function $b$ is called even if $b(x)=b(-x)$.

Our motivation to prove in the following sections identities
for the above mentioned determinants comes
from several problems. The best known problem, called the Fisher-Hartwig 
conjecture, concerns the asymptotics of the determinants of 
Toeplitz matrices for singular symbols. One would like to be able to 
compute the asymptotics of the determinant of $T_N(a)$
when the symbol $a$ has jump 
discontinuities, zeros, or other singularities of a certain form. 
A history of this problem and many known results and applications can
be found in \cite{BS} or \cite{E}. In section five of this paper we 
prove some Fisher-Hartwig type results for certain skewsymmetric
Toeplitz matrices. 

Another interesting problem is to compute asymptotically the
determinants of the matrices $$T_N(a) + H_N(a)$$ where the symbol $a$
also has singularities.  The interest in these asymptotics, especially
in the case where $a$ is even, arose in random matrix theory (see
\cite{BE1} and the references therein). The determination
of these asymptotics will be done in a forthcoming paper \cite{BE2}.

Finally, Hankel matrices defined by the moments of a
function given on a line segment play an important role in orthogonal
polynomial theory and again in
random matrix theory. We refer the reader to \cite{Sz} for orthogonal
polynomial connections and to \cite{M} for a general account of random
matrix theory. In section five we prove two results for the asymptotics
of the determinants of the Hankel moment matrices. These results
allow the function $b$ to have jump discontinuities, but require
the function to be even.

The paper is organized as follows. Sections 2, 3, and 4 contain
all the linear algebra type results which prove the exact identities
for the various types of matrices and are self-contained. The 
asymptotic results are contained in section 5 and use the results
of the previous sections and some already known results for
Toeplitz matrices.


\section{Hankel determinants versus symmetric Toeplitz + Hankel
determinants}

We begin with a preliminary result which will allow us  
to show the relationship with symmetric Toeplitz plus Hankel matrices and
the Hankel moment matrices.

\begin{proposition}\label{p1.1}
Let $\{a_n\}_{n=-\iy}^\iy$ be a sequence of complex numbers such that
$a_n=a_{-n}$ and let $\{b_n\}_{n=1}^\iy$ be a sequence defined by
\bqn\label{f.b}
b_n &=&
\sum_{k=0}^{n-1} {n-1\choose k}(a_{1-n+2k}+a_{2-n+2k}).
\eqn
Define the one-sided infinite matrices
\ben\label{f.AB}
A\;=\;\left(a_{j-k}+a_{j+k+1}\right)_{j,k=0}^\iy,\qquad
B\;=\;\left(b_{j+k+1}\right)_{j,k=0}^\iy,
\een
and the upper triangular one-sided infinite matrix
\bqn\label{f.D}
D &=&
\left(\ba{cccc}
\xi(0,0) & \xi(1,1) & \xi(2,2) & \dots \\
& \xi(1,0) & \xi(2,1) & \dots \\
&& \xi(2,0) & \\
0 &&& \ddots \ea\right)
\quad\mbox{ where }\quad
\xi(n,k)\;=\;{n\choose [\frac{k}{2}]}.
\eqn
Then $B=D^TAD$.
\end{proposition}
\begin{proof}
The assertion is equivalent to the statement that for all
$n,m\ge0$ the following identity holds:
\bqn\label{f1.b1}
b_{n+m+1} &=& \sum_{j=0}^{n}\sum_{k=0}^{m}
(a_{n-j-m+k}+a_{n-j+m-k+1})\xi(n,j)\xi(m,k),
\eqn
where $b_{n+m+1}$ is given by 
\bqn\label{f1.b2}
b_{n+m+1} &=& 
\sum_{r=0}^{n+m}{n+m\choose r}(a_{2r-n-m}+a_{2r-n-m+1}).
\eqn
In order to prove this identity it is sufficient to prove that for each
$s\geq0$ the terms $a_{s}=a_{-s}$ occur as many times in (\ref{f1.b1}) as in
(\ref{f1.b2}). In fact, $a_s$ and $a_{-s}$ occurs in (\ref{f1.b1}) exactly
$N_1+N_2+N_3$ times if $s\ge1$ and $N_1=N_2$ times if $s=0$, where
\bqn
N_1 &=&  \sum_{\sumarr{0\le j\le n}{0\le k\le m}{j-k=n-m-s}}
{n\choose [\frac{j}{2}]}{m\choose [\frac{k}{2}]}
\;\;=\;\;\sum_{\sumarr{0\le j\le n}{m+1\le k\le 2m+1}{j+k=n+m-s+1}}
{n\choose [\frac{j}{2}]}{m\choose [\frac{k}{2}]},
\nn\\
N_2 &=&  \sum_{\sumarr{0\le j\le n}{0\le k\le m}{j-k=n-m+s}}
{n\choose [\frac{j}{2}]}{m\choose [\frac{k}{2}]}
\;\;=\;\;\sum_{\sumarr{n+1\le j\le 2n+1}{0\le k\le m}{j+k=n+m-s+1}}
{n\choose [\frac{j}{2}]}{m\choose [\frac{k}{2}]},
\nn\\
N_3 &=& \sum_{\sumarr{0\le j\le n}{0\le k\le m}{j+k=n+m+1-s}}
{n\choose [\frac{j}{2}]}{m\choose [\frac{k}{2}]}.\nn
\eqn
In the expression for $N_1$ we have made a change of variables
$k\mapsto 2m+1-k$ and in $N_2$ a change of variables $j\mapsto 2n+1-j$.
Hence it follows that
\bqn
N_1+N_2+N_3 &=& 
\sum_{{j,k\ge0 \atop j+k=n+m+1-s}}
{n\choose [\frac{j}{2}]}{m\choose [\frac{k}{2}]}.\nn
\eqn
Moreover, $N_1=N_2=\frac{N_1+N_2+N_3}{2}$ for $s=0$ since then $N_3=0$.

On the other hand, $a_s$ and $a_{-s}$ occurs in (\ref{f1.b2}) exactly 
$M_1+M_2$ times if $s\ge1$ and $M_1=M_2$ times if $s=0$, where
$$
M_1\;\;=\;\; {n+m\choose[\frac{n+m+s}{2}]},\qquad
M_2\;\;=\;\; {n+m\choose[\frac{n+m-s}{2}]}.
$$
Thus we are done as soon as we have shown that
$M_1+M_2=N_1+N_2+N_3$ for each $s\ge0$.

We distinguish two cases. If $n+m+1-s$ is even, then we substitute
$j\mapsto2j$, $k\mapsto 2k$, and $j\mapsto2j+1$, $k\mapsto2k+1$ in
the above expression for $N_1+N_2+N_3$ and arrive at
\bqn
N_1+N_2+N_3 &=&
\sum_{{j,k\ge0 \atop 2j+2k=n+m+1-s}}
{n\choose j}{m\choose k} +
\sum_{{j,k\ge0 \atop 2j+2k=n+m-1-s}}
{n\choose j}{m\choose k}.\nn\\
&=& {n+m\choose\frac{n+m+1-s}{2}}+
{n+m\choose\frac{n+m-1-s}{2}}
\;\;=\;\;M_1+M_2.\nn
\eqn
If $n+m+1-s$ is odd, then we substitute
$j\mapsto2j$, $k\mapsto 2k+1$, and $j\mapsto2j+1$, $k\mapsto2k$ in
the expression for $N_1+N_2+N_3$ and obtain
\bqn
N_1+N_2+N_3 &=& 
2\sum_{{j,k\ge0 \atop 2j+2k=n+m-s}}
{n\choose j}{m\choose k}\nn\\
&=& 2{n+m\choose\frac{n+m-s}{2}}
\;\;=\;\;M_1+M_2,\nn
\eqn
which also completes the proof.
\end{proof}

\begin{theorem}\label{t1.2}
Let $\{a_n\}_{n=-\iy}^\iy$ and $\{b_n\}_{n=1}^\iy$
fulfill the assumptions of the previous proposition. For $N\ge1$
define the matrices
\ben\label{f.ABN}
A_N\;=\;\left(a_{j-k}+a_{j+k+1}\right)_{j,k=0}^{N-1},\qquad
B_N\;=\;\left(b_{j+k+1}\right)_{j,k=0}^{N-1}.
\een
Then $\det A_N=\det B_N$.
\end{theorem}
\begin{proof}
$A_N$ and $B_N$ are the $N\times N$ sections of the infinite
matrices $A$ and $B$ of the previous proposition. Let
$D_N$ be the $N\times N$ sections of the infinite matrix $D$.
Because of the triangular structure of $D$, it follows that
$B_N=D_N^TA_ND_N$. Noting that the entries on the diagonal of
$D$ are equal to $\xi(n,0)=1$, we obtain the desired
assertion.
\end{proof}

The previous theorem shows the connection between the determinants of a 
symmetric
Toeplitz + Hankel matrix on the one hand and a Hankel determinant
on the other hand. We now express
this relationship by using the standard notation for these matrices.

\begin{theorem}\label{t1.3}
Let $a\in L^1(\T)$ be an even function, and define $b\in L^1[-1,1]$ by
\bqn\label{f1.10}
b(\cos\theta) &=& a(e^{i\theta})
\sqrt{\frac{1+\cos\theta}{1-\cos\theta}}.
\eqn
Then $\det(T_N(a)+H_N(a))=\det H_N[b]$.
\end{theorem}
\begin{proof}
The moments of $b$ are given by
\bqn
b_n &=&
\frac{1}{\pi}\int_{-1}^1b(x)(2x)^{n-1}\,dx\nn\\
&=&
\frac{1}{\pi}\int_{0}^\pi a(e^{i\theta})(1+\cos\theta)(2\cos\theta)^{n-1}\,
d\theta\nn\\
&=&
\frac{1}{2\pi}\int_{-\pi}^{\pi}
a(e^{i\theta})(1+e^{-i\theta})(e^{i\theta}+e^{-i\theta})^{n-1}\,
d\theta\nn\\
&=&
\frac{1}{2\pi}\int_{-\pi}^{\pi}a(e^{i\theta})
\left(\sum_{k=0}^{n-1}(e^{i(n-1-2k)\theta}+e^{i(n-2-2k)\theta})
{n-1\choose k}\right)\,d\theta\nn\\
&=&
\sum_{k=0}^{n-1}{n-1\choose k}\left(a_{-n+1+2k}+a_{-n+2+2k}\right).\nn
\eqn
Here we have made a change of variables $x=\cos\theta$ and written
$(e^{i\theta}+e^{-i\theta})^{n-1}$ using the binomial formula.
With regard to (\ref{f.b}) and Theorem \ref{t1.2} this
completes the proof.
\end{proof}

In regard to relation (\ref{f1.10}) we remark that $b\in L^1[-1,1]$
in and only if $a(e^{i\theta})(1+\cos\theta)\in L^1(\T)$.

Thus at this point we have shown that if $a$ and $b$ satisfy the relation
(2.12), then
$$\det H_N[b] = \det(T_N(a)+H_N(a)).$$ But actually more can be done in the case 
that the symbol $a$ satisfies a quarter wave symmetry property. Then, in fact, 
certain Hankel moment determinants can be written as Toeplitz determinants.
The symbol $b(x)\in L^1[-1,1]$ of these Hankel determinants is of the form
\bqn\label{f1.11}
b(x)&=& b_0(x)\sqrt{\frac{1+x}{1-x}}
\eqn
where $b_0(-x)=b_0(x)$ for all $x\in[-1,1]$.

We first begin with the following auxiliary result.
In what follows, let $W_N$ stand for the matrix acting on $\C^N$ by
$$
W_N:(x_0,x_1,\dots,x_{N-1})\mapsto (x_{N-1},\dots,x_1,x_0),
$$
and let $I_N$ denote the $N\times N$ identity matrix.

\begin{proposition}\label{p1.4new}
Let $a\in L^1(\T)$ and assume that $a(-t)=a(t\iv)=a(t)$. Define
\bqn
d(e^{i\theta})&=& a(e^{i\theta/2}).
\eqn
Then $\det(T_{N}(a)+H_N(a))=\det T_N(d)$.
\end{proposition}
\begin{proof}
Note first that $d(t)$ is well defined since $a(t)=a(-t)$. Moreover,
$a_{2n+1}=0$ and $a_{2n}=d_n$. By rearranging rows and columns of
$T_N(a)+H_N(a)$ in an obvious way, it is easily seen that this matrix
is similar to
\bqn
\left(\ba{cc}\left(a_{2j-2k}\right)_{j,k=0}^{N_1-1}&0\\ 0&
\left(a_{2j-2k}\right)_{j,k=0}^{N_2-1}\ea\right)+
\left(\ba{cc}0&\left(a_{2j+2k+2}\right)_{j=0,\hspace{1.4ex}k=0}^{N_1-1,N_2-1}\\
\left(a_{2j+2k+2}\right)_{j=0,\hspace{1.4ex}k=0}^{N_2-1,N_1-1}&0\ea\right)\nn
\eqn
where $N_1=\left[\frac{N+1}{2}\right]$ and $N_2=\left[\frac{N-1}{2}\right]$.
This matrix equals
$$\left(\ba{cc} T_{N_1}(d)&H_{N_1,N_2}(d)\\H_{N_2,N_1}(d)&T_{N_2}(d)\ea
\right),$$
where $H_{N_1,N_2}(d)$ and $H_{N_2,N_1}(d)$ are Hankel matrices of size
$N_1\times N_2$ and $N_2\times N_1$, respectively.
Multiplying the last matrix from the left and the right with the
diagonal matrix $\diag(W_{N_1},I_{N_2})$ we obtain the matrix
$T_N(d)$. Notice in this connection that $d_{n}=d_{-n}$
since $a(t\iv)=a(t)$.
\end{proof}

\begin{corollary}\label{c1.5}
Let $b\in L^1[-1,1]$ and suppose (\ref{f1.11}) holds with
$b_0(-x)=b_0(x)$ for all $x\in[-1,1]$.
Define the function
\bqn
d(e^{i\theta})&=&b_0(\cos\frac{\theta}{2}).
\eqn
Then $\det H_N[b]=\det T_{N}(d)$.
\end{corollary}
\begin{proof}
Since $b_0(x)=b_0(-x)$ it follows from definition 
(\ref{f1.10}) that $a(-t)=a(t\iv)=a(t)$. Now we can apply 
Theorem \ref{t1.3} and Proposition \ref{p1.4new} in order to obtain the
identity $\det H_N[b]=\det (T_N(a)+H_N(a))=\det T_N(d)$.
\end{proof}
Concerning the previous corollary, we wish to emphasize that the function
$d$ is even, and hence the matrices $T_N(d)$ are symmetric.


\section{Symmetric Toeplitz + Hankel determinants versus
skewsymmetric Toeplitz determinants}

The main result of this section has been established in \cite[Lemma 18]{Kra}
and proved in \cite[Lemma 1]{Gor} and \cite[Proof of Thm.~7.1(a)]{Stem}.
We give a slightly simplified and self-contained proof here.

\begin{theorem}\label{t2.3}
Let $\{a_n\}_{n=-\iy}^{\iy}$ be a sequence of complex numbers such that
$a_{-n}=a_{n}$. Let $c_{n}$ be defined by 
\bqn\label{f.cn}
c_{n} &=& \sum_{k=-n+1}^n a_{k}\quad\mbox{ for }
n>0,
\eqn
and put $c_0=0$ and $c_{-n}=-c_{n}$. 
Then $\det T_{2N}(c)=(\det(T_N(a)+H_N(a)))^2$.
\end{theorem}
\begin{proof}
First of all we
multiply the matrix $T_{2N}(c)$ from the left and right with
$\diag(W_N,I_N)$. We obtain the matrix
$$
\left(\ba{cc} T_N(\tilde{c})&H_N(\tilde{c})\\H_N(c)&T_N(c)\ea\right)
\;\;=\;\;
\left(\ba{cc}-T_N(c)&-H_N(c)\\H_N(c)&T_N(c)\ea\right)
$$
by observing that $\tilde{c}=-c$. Next we claim that
\bqn
&&\hspace*{-5ex}
\left(\ba{cc}T_N(1-t)&0\\T_N(t)&I_N\ea\right)
\left(\ba{cc}-T_N(c)&-H_N(c)\\H_N(c)&T_N(c)\ea\right)
\left(\ba{cc}T_N(1-t\iv)&T_N(t\iv)\\0&I_N\ea\right)\nn\\
&=&
\left(\ba{cc}I_N&0\\0&T_N(1+t)\ea\right)
\left(\ba{cc}X_N&-T_N(a)-H_N(a)\\T_N(a)+H_N(a)&0\ea\right)
\left(\ba{cc}I_N&0\\0&T_N(1+t\iv)\ea\right)\nn
\eqn
with a certain matrix $X_N$. If we take the determinant of this equation,
we obtain the desired determinant identity.

In order to proof the above matrix identity it suffices to show
that the following three equations hold:
\bqn
T_N(c)-T_N(t)T_N(c)T_N(t\iv)+H_N(c)T_N(t\iv)-T_N(t)H_N(c)
\;\;=\;\;0,\hspace*{-24ex}&&\label{f2.I}\\[1ex]
-T_N(t)T_N(c)T_N(1-t\iv)+H_N(c)T_N(1-t\iv) &=&
T_N(1+t)\left(T_N(a)+H_N(a)\right),\label{f2.II}\\[1ex]
T_N(1-t)T_N(c)T_N(t\iv)+T_N(1-t)H_N(c)&=&
\left(T_N(a)+H_N(a)\right)T_N(1+t\iv).\qquad\label{f2.III}
\eqn
Notice that (\ref{f2.III}) can be obtained from (\ref{f2.II}) by passing to
the transpose. Moreover, by employing (\ref{f2.I}) equation (\ref{f2.II})
reduces to
\bqn
T_N(1-t)\left(T_N(c)+H_N(c)\right) &=&
T_N(1+t)\left(T_N(a)+H_N(a)\right).\label{f2.IV}
\eqn

Let us first prove (\ref{f2.I}). We introduce the $N\times1$ column vectors
$e_0=(1,0,0,\dots,0)^T$ and $\gamma_N=(0,c_1,c_2,\dots,c_{N-1})^T$.
Then
$$
T_N(c)-T_N(t)T_N(c)T_N(t\iv) \;\;=\;\; \gamma_Ne_0     ^T-e_0\gamma_N^T
\;\;=\;\; T_N(t)H_N(c)-H_N(c)T_N(t\iv),
$$
whence indeed (\ref{f2.I}) follows. 

Next we remark that from the definition of the
sequences $\{a_n\}_{n=-\iy}^\iy$ and $\{c_n\}_{n=-\iy}^\iy$
it follows that $c_n-c_{n-1}=a_n+a_{n-1}$ for all $n\in\Z$.
Introducing the column vectors
$\hat{\gamma}_N=(c_1,\dots,c_N)^T$, $\alpha_N=(a_0,\dots,a_{N-1})^T$ and
$\hat{\alpha}_{N}=(a_1,\dots,a_N)^T$, it can be readily verified that
\bqn
T_N(1-t)T_N(c) &=& \left(c_{j-k}-c_{j-k-1}\right)_{j,k=0}^{N-1}
-e_0\hat{\gamma}_N^T,\nn\\
T_N(1+t)T_N(a) &=& \left(a_{j-k}+a_{j-k-1}\right)_{j,k=0}^{N-1}
-e_0\hat{\alpha}_N^T,\nn\\
T_N(1-t)H_N(c) &=& \left(c_{j+k+1}-c_{j+k}\right)_{j,k=0}^{N-1}
+e_0\gamma_N^T,\nn\\
T_N(1+t)H_N(a) &=& \left(a_{j+k+1}+a_{j+k}\right)_{j,k=0}^{N-1}
-e_0\alpha_N^T.\nn
\eqn
Using the above relation $c_n-c_{n-1}=a_n+a_{n-1}$, it follows that
\bqn
T_N(1-t)T_N(c)-T_N(1+t)T_N(a)
&=& -e_0\hat{\gamma}^T_N+e_0\hat{\alpha}^T_N\nn\\
T_N(1+t)H_N(a)-T_N(1-t)H_N(c)
&=&-e_0\alpha^T_N-e_0\gamma^T_N\nn.
\eqn
Since $\hat{\gamma}_N-\gamma_N=\hat{\alpha}_N+\alpha_N$ by the same
relation, this implies equation (\ref{f2.IV}).
\end{proof}

The results of this theorem are not easy to rephrase by using the
classical notation for Toeplitz and Hankel matrices. Consider, for instance,
the simplest case where $a(t)\equiv1$. Then $c_n={\rm sign}(n)$
which are not the Fourier coefficients of an $L^1$-function.
For more information on how one can nevertheless express the relationship
between the symbols $a$ and $c$, and how the asymptotics for certain
of the above determinants can be determined we refer to \cite{BE2}.


\section{Hankel determinants versus skewsymmetric Toeplitz determinants}

The results of the previous two sections allow us to
establish an identity between Hankel determinants and determinants of
skewsymmetric Toeplitz matrices. The next theorem is an additional needed
ingredient for the identity.

\begin{theorem}
Let $\{c_n\}_{n=-\iy}^{\iy}$ be a sequence of complex numbers
such that $c_{-n}=-c_n$ for all $n\in\Z$. Define numbers
$\{b_n\}_{n=1}^\iy$ by
\bqn\label{f.bn-cn}
b_n &=& \sum_{k=0}^{\left[\frac{n}{2}\right]}
\left\{{n-1 \choose k}-{n-1 \choose k-1}\right\}c_{n-2k}.
\eqn
Moreover, define the matrices
$$
B_N=\left(b_{j+k+1}\right)_{j,k=0}^{N-1},\qquad
C_{2N}=\left(c_{j-k}\right)_{j,k=0}^{2N-1}.
$$
Then $\det C_{2N}=(\det B_N)^2$.
\end{theorem}
\begin{proof}
In formula (\ref{f.cn}) the numbers $c_n$ are defined in terms of
the numbers $a_{-n+1},\dots,a_{n}$. By a simple inspection of this formula,
it is easy to see that for any given sequence $\{c_n\}_{n=-\iy}^\iy$ there
exists a sequence $\{a_n\}_{n=-\iy}^\iy$ such that (\ref{f.cn})
and $a_{n}=a_{-n}$ holds for all positive $n$.

Now let us define the numbers $b_n$ not by (\ref{f.bn-cn}) but by
(\ref{f.b}). Then with $B_{N}$ and $C_{2N}$ defined as above
it follows from Theorem \ref{t1.2} and Theorem \ref{t2.3} that
$\det C_{2N}=(\det B_N)^2$. It remains to show that (\ref{f.bn-cn})
holds.

Indeed, we have that
\bqn
\lefteqn{\sum_{k=0}^{\left[\frac{n}{2}\right]}
\left\{{ n-1 \choose k}-{n-1 \choose k-1}\right\} c_{n-2k}
\;\;=\;\;
\sum_{k=0}^{\left[\frac{n}{2}\right]}
\left\{{n-1 \choose k}-{n-1 \choose k-1}\right\}
\sum_{j=-n+2k+1}^{n-2k} a_j}
\hspace*{15ex}\nn\\
&=&
\sum_{-n+2k+1\le j\le n-2k\atop0\le 2k\le n}
\left\{{n-1 \choose k}-{n-1 \choose k-1}\right\}a_j\nn\\
&=&
\sum_{j=-n+1}^{n}
\sum_{k=0}^{\min\left\{\left[\frac{n-j}{2}\right],
\left[\frac{n+j-1}{2}\right]\right\}}
\left\{{n-1 \choose k}-{n-1 \choose k-1}\right\}a_j\nn\\
&=&
\sum_{j=-n+1}^{n}
{n-1\choose\min\left\{\left[\frac{n-j}{2}\right],
\left[\frac{n+j-1}{2}\right]\right\}}a_j\nn\\
&=&
\sum_{j=-n+1}^n{n-1\choose\left[\frac{n-j}{2}\right]}a_j
\;\;=\;\;
\sum_{k=0}^{n-1}{n-1\choose k}(a_{2k+1-n}+a_{2k+2-n}).\nn
\eqn
By formula (\ref{f.b}) this is equal to $b_n$.
\end{proof}

We again express the above relationship in
terms of the standard notation.

\begin{theorem}\label{t3.2}
Let $b\in L^1[-1,1]$ and define $c\in L^1(\T)$ by
\bqn\label{f3.15}
c(e^{i\theta})&=& i\,{\rm sign}(\theta)\,b(\cos\theta),\qquad
-\pi<\theta<\pi.
\eqn
Then $\det T_{2N}(c)=(\det H_N[b])^2$.
\end{theorem}
\begin{proof}
Obviously, $c(e^{-i\theta})=-c(e^{i\theta})$. Hence $c_{-n}=-c_{n}$.
It is sufficient to verify formula (\ref{f.bn-cn}) for the 
Fourier coefficients and moments.
First of all,
\bqn
c_n &=&
\frac{1}{\pi}\int_0^{\pi}b(\cos\theta)\sin(n\theta)\,d\theta.
\nn\eqn
Hence
\bqn
b_n &=&
\frac{1}{\pi}\int_0^{\pi}b(\cos\theta)
\left(\sum_{k=0}^{\left[\frac{n}{2}\right]}
\left\{{n-1 \choose k}-{n-1 \choose k-1}\right\}\sin((n-2k)\theta)
\right)d\theta.\nn
\eqn
The expression in the big braces equals (by a change of variables
$k\mapsto n-k$ in the second part of the sum)
\bqn
\lefteqn{
\sum_{k=0}^{\left[\frac{n}{2}\right]}{n-1 \choose k}\sin((n-2k)\theta)
-\sum_{k=n-\left[\frac{n}{2}\right]}^n
{n-1 \choose n-k-1}\sin((2k-n)\theta)
}\hspace{10ex}
\nn\\
&=&
\sum_{k=0}^{n-1} {n-1 \choose k}\sin((n-2k)\theta)
\;\;=\;\;
(2\cos\theta)^{n-1}\sin\theta.\nn
\eqn
Hence
\bqn
b_n &=&
\frac{1}{\pi}\int_0^\pi
b(\cos\theta)(2\cos\theta)^{n-1}\sin\theta\,d\theta.\nn
\eqn
Now it is easy to see that $b_n$ are the moments of the function $b$.
\end{proof}

Regarding relation (\ref{f3.15}) we remark that
$c\in L^1(\T)$ if and only if $b(x)/\sqrt{1-x^2}\in L^1[-1,1]$.

At this point we have three main identities for Hankel moment determinants,
one which follows from Theorem 2.3, one which follows from Corollary
2.5 and finally one which follows from the previous theorem. If we 
desire to find the asymptotics of the determinants of the Hankel moment
matrices it is clear that the corresponding asymptotics for
Toeplitz matrices need to be derived. In particular, in light of
Theorem 4.2 and formula (4.22), it is desirable to 
compute the asymptotics of the
Toeplitz determinant $\det T_{2N}(c)$, where $c$ satisfies
$c(e^{-i\theta})=-c(e^{i\theta})$ and accordingly implies
that the Toeplitz matrices are skewsymmetric. 
Note from this it follows that $\det T_{2N+1}(c)=0$ for all
$N$. However, this implies that a single asymptotic formula for
the determinants, such as the one given in the
classical Szeg\"{o} limit theorem,
or the more general Fisher-Hartwig formulas would not make sense here.
In the following section we nevertheless compute the asymptotics of such
Toeplitz determinants in some cases and raise a conjecture about 
more general cases.


\section{Asymptotics of certain skewsymmetric Toeplitz determinants
and Hankel determinants}

Our goal of this section is to consider Toeplitz determinants with
generating function $c(e^{i\theta}) = \chi(e^{i\theta})a(e^{i\theta})$
where $a$ is an even functions and 
\bqn
\chi(e^{i\theta}) &=& i\,{\rm sign}(\theta),
\qquad -\pi<\theta<\pi.
\eqn
Let $t_{\beta}(e^{i\theta})$ stand for the function
\bqn
t_{\beta}(e^{i\theta}) &=& e^{i\beta(\theta-\pi)},
\qquad 0<\theta<2\pi.
\eqn
This function has a single jump at $t=1$ whose size is determined by
the parameter $\beta$.

In the following proposition we assume that $a$ is not necessarily an even
function but satisfies instead a rotation symmetry condition.

\begin{proposition}\label{p3.3}
Assume that $a\in L^1(\T)$ satisfies the relation $a(-t)=a(t)$
for $t\in\T$. Define the functions
\bqn
d(e^{i\theta})=a(e^{i\theta/2}),\qquad
d_1(e^{i\theta})=t_{-1/2}(e^{i\theta})d(e^{i\theta}),\qquad
d_2(e^{i\theta})=t_{1/2}(e^{i\theta})d(e^{i\theta}).\nn
\eqn
Then $\det T_{2N}(a)=(\det T_N(d))^2$ and
$\det T_{2N}(\chi a)=\det T_{N}(d_1)\det T_{N}(d_2)$.
\end{proposition}
\begin{proof}
{}From the assumptions $a(t)=a(-t)$ it follows that the Fourier coefficients
$a_{2n+1}$ are zero. Hence $T_{2N}(a)$ has a checkered pattern, and
rearranging rows and columns it is easily seen that $T_{2N}(a)$ is similar
to the matrix $\diag(T_N(d),T_N(d))$.

The Fourier coefficients $c_{2n}$ of $c(t)=\chi(t)a(t)$ are equal to zero.
By rearranging the rows and columns of $T_{2N}(\chi a)$ in the same way
as above it becomes apparent that $T_{2N}(\chi a)$ is similar to a matrix
$$
\left(\ba{cc}0&D_2\\D_1&0\ea\right)\quad\mbox{ where }
D_{1}=\left(c_{2(j-k)+1}\right)_{j,k=0}^{N-1}\mbox{ and }
D_{2}=\left(c_{2(j-k)-1}\right)_{j,k=0}^{N-1}.
$$
{}From the identity
\bqn\label{f4.23}
\chi(e^{i\theta}) \;\;=\;\;
t_{-1/2}(e^{i\theta})t_{1/2}(e^{i(\theta-\pi)}) \;\;=\;\;
-t_{1/2}(e^{i\theta})t_{-1/2}(e^{i(\theta-\pi)})
\eqn
it follows that $d_1(e^{i\theta})=e^{-i\theta/2}c(e^{i\theta/2})$ and
$d_2(e^{i\theta})=-e^{i\theta/2}c(e^{i\theta/2})$.
Hence $D_1=T_{N}(d_1)$ and $D_2=-T_N(d_2)$. Since
$\det T_{2N}(c)=(-1)^N \det D_1 \det D_2$, this completes the proof.
\end{proof}

Hence we have reduced the computation of $\det T_{2N}(\chi a)$
to the Toeplitz determinants $T_N(d_1)$ and $T_N(d_2)$,
for which in the case of piecewise continuous functions it is
possible to apply the Fisher-Hartwig conjecture under certain assumptions.

The following result, which is taken from \cite{E}, makes this explicit.
Therein $G(\cdot)$ is the Barnes $G$-function \cite{WW}, 
$d_{0,\pm}$ are the Wiener-Hopf factors of the function $d_0$,
\bqn
d_{0,\pm}(e^{i\theta})&=&\exp\left(\sum_{k=1}^\iy
[\log d_0]_{\pm k}e^{\pm ik\theta}\right),
\eqn
and
\bqn
d_\pm(e^{i\theta})&=& d_{0,\pm}(e^{i\theta})\prod_{r=1}^R
\Big(1-e^{\pm i(\theta-\theta_r)}\Big)^{\pm\beta_r}
\eqn
are the generalized Wiener-Hopf factors of $d$.

\begin{proposition}
Let 
\bqn\label{f4.25}
d(e^{i\theta}) &=& d_0(e^{i\theta})\prod_{r=1}^R
t_{\beta_r}(e^{i(\theta-\theta_r)}),
\eqn
where $d_0$ is an infinitely differentiable nonvanishing function with
winding number zero, $\theta_1,\dots,\theta_R\in(0,2\pi)$ are distinct
numbers, and $\beta_1,\dots,\beta_R$ are complex parameters satisfying
$|\Re\beta_r|<1/2$ for all $r=1,\dots,R$. Then
\bqn\label{f3.26}
\frac{\det T_N(t_{-1/2}d)}{\det T_N(d)} &\sim&
N^{-1/4}G(1/2)G(3/2)d_+(1)^{-1/2}d_-(1)^{1/2},
\qquad N\to\iy,\nn\\
\frac{\det T_N(t_{1/2}d)}{\det T_N(d)} &\sim&
N^{-1/4}G(1/2)G(3/2)d_+(1)^{1/2}d_-(1)^{-1/2},
\qquad N\to\iy.\nn
\eqn
Moreover,
\bqn
\det T_N(d) &\sim& 
\G^NN^{\Omega}E,\qquad N\to\iy,
\eqn
where $\G=\exp\left(\frac{1}{2\pi}\int_{0}^{2\pi}
\log d_0(e^{i\theta})\,d\theta\right)$,
$\Omega=-\sum\limits_{r=1}^R\beta_r^2$, and $E$ is another constant.
\end{proposition}
(The constant $E$ is quite complicated, so in the interest of
brevity, we omit the exact formula from this paper and refer
to \cite{BT,BS,E} for an explicit representation.)

The previous propositions yield the following results.
We keep the same notation.

\begin{corollary}\label{c4.3}
Let $d$ be a function of the form (\ref{f4.25}) and assume that the
same conditions as above are fulfilled. Let $a(e^{i\theta})=d(e^{2i\theta})$.
Then
\bqn\label{f4.29}
\frac{\det T_{2N}(\chi a)}{\det T_{2N}(a)} &\sim&
N^{-1/2}G^2(1/2)G^2(3/2),\qquad N\to\iy,
\eqn
and
\bqn
\det T_{2N}(a) &\sim& \G^{2N}N^{2\Omega}E^2,\qquad N\to\iy.
\eqn
\end{corollary}

The following corollary gives an asymptotic formula for the determinants of
Hankel moment matrices in the special case where the symbol is even.
\begin{corollary}\label{c4.4}
Let $b \in L^1[-1,1]$ such that $b(-x)= b(x).$ Define
$d(e^{i\theta}) = b(\cos(\theta/2))$ and suppose that $d$ is of
the form (\ref{f4.25}). 
Then
\bqn
\det H_N[b] &\sim& \G^N N^{\Omega -1/4}G(1/2)G(3/2)E,\qquad N\to\iy.
\eqn
\end{corollary}
\begin{proof}
Define $a(e^{i\theta})=b(\cos\theta)$. Then Theorem \ref{t3.2} implies
that $(\det H_N[b])^2=\det T_N(\chi a)$. Since $b(x)=b(-x)$ the
function is well defined and $a(e^{i\theta})=d(e^{2i\theta})$.
Now the formula follows from  Corollary \ref{c4.3} and by
taking square roots.
\end{proof}

The interesting point in  Corollary \ref{c4.3} is that the asymptotic
limit of (\ref{f4.29}) does not depend on the underlying function $a$.
We remark that we have proved this limit relation for certain 
piecewise continuous functions $a$ subject to the condition
$a(-t)=a(t)$. Our primary goal was however to determine the limit for
certain functions $a$ satisfying the relation $a(t\iv)=a(t)$.
Our conjecture is that the asymptotic limit is given by the above 
expression in general also for those functions.

In order to support this hypothesis we resort to the generalization of the
Fisher-Hartwig conjecture, which has not yet been proved, but is strongly
suggested by examples. Since $\det T_{2N+1}(\chi a)=0$
for all $N$ (under the assumption $a(t\iv)=a(t)$), the asymptotics 
of $T_N(\chi a)$ can only be described by the generalized but
not the original conjecture. The crucial observation is that one has
several possibilities for representing $\chi a$ in a form like
(\ref{f4.25}). Indeed, from (\ref{f4.23}) it follows
that
$$
\chi(e^{i\theta}) a(e^{i\theta}) \;\;=\;\;
t_{-1/2}(e^{i\theta})t_{1/2}(e^{i(\theta-\pi)})a(e^{i\theta})
\;\;=\;\;
-t_{1/2}(e^{i\theta})t_{-1/2}(e^{i(\theta-\pi)})a(e^{i\theta}),
$$
tacitly assuming that $a$ admits also representation of the form
(\ref{f4.25}) with appropriate properties.

Then the generalized conjecture predicts \cite{BT,E} that
\bqn
\det T_{N}(\chi a) &\sim&
\det T_N(t_{-1/2}(e^{i\theta}))\det T_N(t_{1/2}(e^{i(\theta-\pi)}))
\det T_N(a) E_1\nn\\&&
+(-1)^N\det T_N(t_{1/2}(e^{i\theta}))\det T_N(t_{-1/2}(e^{i(\theta-\pi)}))
\det T_N(a) E_2,\nn
\eqn
where $E_1$ and $E_2$ are the ``correlation'' constants
\bqn
E_1 &=& E(t_{-1/2}(e^{i\theta}),t_{1/2}(e^{i(\theta-\pi)}))
E(t_{-1/2}(e^{i\theta}),a)E(t_{1/2}(e^{i(\theta-\pi)}),a)
\nn\\
E_2 &=& E(t_{1/2}(e^{i\theta}),t_{-1/2}(e^{i(\theta-\pi)}))
E(t_{1/2}(e^{i\theta}),a)E(t_{-1/2}(e^{i(\theta-\pi)}),a)
\nn
\eqn
with $E(\cdot,\cdot)$ defined by
\bqn
E(b,c) &=& \exp\left(\lim_{r\to1-0}\sum_{k=1}^\iy\Big(
k[\log h_rb_+]_k[\log h_rc_-]_{-k}+
k[\log h_rb_-]_{-k}[\log h_rc_+]_k\Big)\right),\nn
\eqn
$h_rb_\pm$ and $h_rc_\pm$ denoting the harmonic extensions of the
Wiener-Hopf factors of $b_\pm$ and $c_\pm$.

{}From all this it follows that
\bqn
\frac{\det T_{2N}(\chi a)}{\det T_{2N}(a)} &\sim&
(2N)^{-1/2}G^2(1/2)G^2(3/2)(E_1+E_2),\nn
\eqn
where a straightforward computation of the constants gives
\bqn
E_1 &=& 2^{-1/2}\left(\frac{a_+(-1)a_-(1)}{a_-(-1)a_+(1)}\right)^{1/2},\nn\\
E_2 &=& 2^{-1/2}\left(\frac{a_+(-1)a_-(1)}{a_-(-1)a_+(1)}\right)^{-1/2}.\nn
\eqn
The assumption that $a(t\iv)=a(t)$ implies that 
$a_-(t)=\gamma a_+(t\iv)$ with a certain constant $\gamma\neq0$.
Hence
\bqn
E_1=E_2=2^{-1/2},\nn
\eqn
which leads to the conjecture that
\bqn
\frac{\det T_{2N}(\chi a)}{\det T_{2N}(a)} &\sim&
N^{-1/2}G^2(1/2)G^2(3/2),\qquad N\to\iy.
\eqn
Using Theorem \ref{t3.2} we arrive at a conjecture for the
Hankel moment matrices:
\bqn
\frac{\det H_{N}[b]}{\sqrt{\det T_{2N}(a)}} &\sim&
N^{-1/4}G(1/2)G(3/2),\qquad N\to\iy,
\eqn
where $a(e^{i\theta})=b(\cos\theta)$. We remark that this formula is in 
accordance with Corollary \ref{c4.4}.

We end this section by noting one other result that
follows from our identities and Corollary \ref{c1.5}. This result applies to
Hankel moment matrices with a special case of Jacobi
weights and computes the asymptotics
for $\det H_N[b]$ where $b$ is of the form $b_0(x)\sqrt{\frac{1+x}{1-x}}$
with an even function $b_0$.

\begin{corollary}
Suppose $b \in L^1[-1,1]$ is of the above form with
an even function $b_0$. Let $d(e^{i\theta})= b_0(\cos(\theta/2))$ and
suppose the $d$ is of the form (\ref{f4.25}).
Then 
$$\det H_N[b] \sim \G^N N^{\Omega}E ,\qquad N\to \iy.$$
\end{corollary}


\end{document}